\documentclass[11pt]{amsart}
\usepackage{amsmath,amsthm,amscd,amssymb, color}
\usepackage{latexsym}
\usepackage[colorlinks, citecolor = blue,backref]{hyperref}
\usepackage[alphabetic]{amsrefs}
\usepackage[capitalize, nameinlink]{cleveref}
\usepackage{enumerate}
\usepackage[margin=1.333in]{geometry}
\usepackage{graphicx}
\usepackage[font=footnotesize]{caption}

\usepackage{array}   
\newcolumntype{L}{>{$}l<{$}} 
\newcolumntype{R}{>{$}r<{$}} 
\newcolumntype{C}{>{$}c<{$}}

\newcommand{\Q}{\mathbb Q}

\newcommand{\Z}{\mathbb Z}
\newcommand{\F}{\mathbb F}

\renewcommand{\phi}{\varphi}

\newcommand{\calO}{\mathcal O}

\newcommand{\bmx}{\left( \begin{matrix}}
\newcommand{\emx}{\end{matrix} \right)}

\newcommand{\new}{\mathrm{new}}

\renewcommand{\mod}{\bmod}

\usepackage{bbm}

\newcommand{\leg}{\overwithdelims ()}

\DeclareMathOperator{\GL}{GL}

\DeclareMathOperator{\tr}{tr}

\crefname{question}{Question}{Questions}

\crefname{conja}{Conjecture}{Conjectures}

\theoremstyle{remark}

\theoremstyle{definition}

\numberwithin{equation}{section}

\begin{document}

\title{Variations on murmurations}

\author{Kimball Martin}
\address{Department of Mathematics, University of Oklahoma, Norman, OK 73019 USA}
\email{kimball.martin@ou.edu}
\address{Department of Mathematics $\cdot$ International Research and Education Center, Graduate School of Science, Osaka Metropolitan University, Osaka 558-8585, Japan}
\email{kimball@omu.ac.jp}

\date{\today}

\maketitle

\begin{abstract} We explore several variations on the recently discovered phenomena of murmurations for elliptic curves and modular forms.
\end{abstract}

\setcounter{tocdepth}{1}
\tableofcontents


\section{Introduction}


Murmurations are a subtle correlation in a sufficiently large family of objects, such as elliptic curves or modular forms, between the root number/rank of the objects and Fourier/Dirichlet coefficients.  They were first discovered numerically in the context of elliptic curves by He, Lee, Oliver, and Pozdnyakov \cite{HLOP}.  Later, extensive computations of Sutherland with ideas of others elucidated the patterns in murmurations, and indicated such murmurations exist for other families of objects as well, such as modular forms and genus 2 curves -- e.g., see slides or data on Sutherland's website.\footnote{\url{https://math.mit.edu/~drew/}.}  Zubrilina has proved the existence and some properties of murmurations in the context of modular forms of squarefree level \cite{zubrilina}.  Murmurations have also been exhibited in several other settings, such as Dirichlet characters \cite{LOP}, Maass forms \cite{maass}, and higher rank zeta functions of elliptic curves \cite{shi-weng}.

In this note, I will explore some different kinds of variants of murmurations, primarily from an empircal perspective.  We will go beyond the paradigm of looking at correlations between root numbers and Fourier/Dirichlet coefficients, and think about the following perspectives:

\begin{enumerate}
\item expected minus actual solution counts
\item averaging arithmetic functions over subsets of integers
\item traces of ``arithmetic'' linear operators, and generalizations
\end{enumerate}

Note that $a_p$'s of elliptic curves are of type (1).  Averages of Fourier coefficients of modular forms are traces of Hecke operators, i.e., type (3), which via the trace formula can be expressed as class number sums, i.e., type (2).  The ``and generalizations'' in (3) refers to considering alternate types of trace formulas, such as Kuznetsov or relative trace formulas, as opposed to Eichler--Selberg-type trace formulas.


\subsection*{Acknowledgements}
This note is largely based on presentations at the workshops
\emph{Murmurations in Arithmetic} (ICERM, July 6--8, 2023) and 
\emph{Murmurations in Arithmetic Geometry and Related Topics} (Simons Center, Nov 11--15, 2024), as well as ensuing discussions with participants.  For this I am grateful to the organizers.  I also thank Alex Cowan, Thomas Oliver, Andrew Sutherland, and Nina Zubrilina for several helpful discussions and comments.

\section{Review of murmurations} \label{sec:review}

Let $\mathcal F$ be one of the following two types of families:

\begin{itemize}
\item $\mathcal E$ - elliptic curves up to isogeny, partially ordered by conductor $N$

\item $\mathcal H_k$ - weight $k$ cuspidal newforms with trivial nebentypus, partially ordered by level $N$
\end{itemize}

Fix $\beta > 1$; by default we take $\beta = 2$ following Sutherland.  Let $\mathcal F^{\pm}(N)$ be the set objects $f \in \mathcal F$ with conductor/level $N$ and root number $\pm 1$.  If $\mathcal F = \mathcal E$, set $k=2$.  Consider the averages
\[ \mathcal A_{\mathcal F}^{\pm}(p,X) = \mathcal A_{\mathcal F}^{\pm}(p,X; \beta) = p^{1 - \frac k2}
\frac{\sum'_{X \le N \le \beta X} \sum_{f \in \mathcal F^{\pm}(N)} a_p(f)}{\sum'_{X \le N \le \beta X} \# \mathcal F^{\pm}(N)}. \]
The $p^{1 - \frac k2}$ is a normalization factor, which is 1 if $\mathcal F = \mathcal E$ or $\mathcal H_2$.
Here the prime on the sums over $N$ denotes a possible restriction on the $N$ considered --- in this section we will take $N$ squarefree and coprime to $p$, largely as a computational convenience.   

\begin{figure}
\begin{minipage}{.5\textwidth}
\captionof{figure}{Murmurations for weight 2 modular forms of squarefree level $1000 \le N \le 2000$}
 \includegraphics[width=\textwidth]{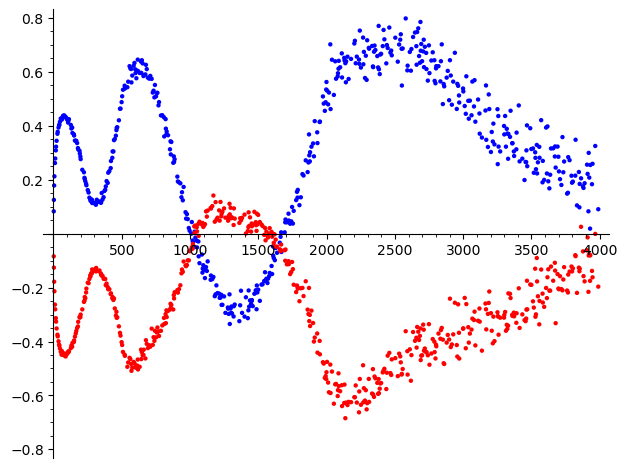}
\label{fig:fixedroot1}
\end{minipage}%
\begin{minipage}{.5\textwidth}
\captionof{figure}{Murmurations for weight 2 modular forms of squarefree levels $2000 \le N \le 4000$}
 \includegraphics[width=\textwidth]{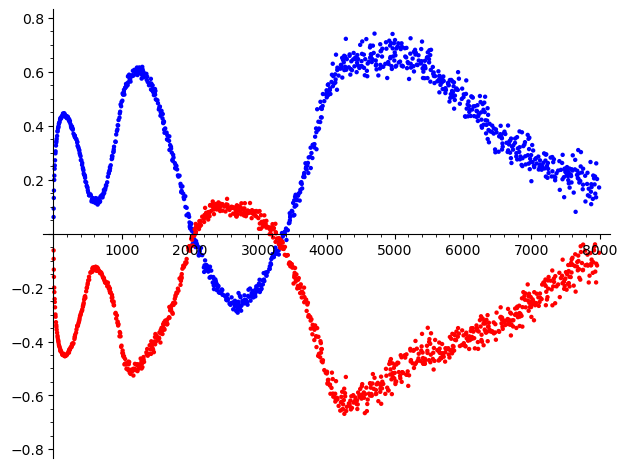}
\label{fig:fixedroot2}
\end{minipage}
\end{figure}

For $\mathcal F = \mathcal H_2$, we plot $\mathcal A_{\mathcal F}^{\pm}(p; X)$ as a function of $p$ for $X=1000$ in \cref{fig:fixedroot1} and $X=2000$ in \cref{fig:fixedroot2}; the case of root number $+1$ is plotted in blue and root number $-1$ in red.  Note that these graphs look essentially the same, even though the underlying data comes from completely disjoint sets of newforms!  We have chosen the horizontal scale to range up to $p \le 4X$ in both cases, so one might conjecture that if we continue to make such graphs for larger and larger $X$, they will tend to the graph of a smooth function in $\frac pX$, called the \emph{murmuration function}.  This limiting property is called \emph{scale invariance} in $\frac pX$ (or $\frac pN$).

Indeed this conjecture holds for any even $k \ge 2$ and any $\beta > 1$, as shown by Zubrilina \cite{zubrilina} with careful estimates of class number sums coming from trace formulas.  Zubrilina also describes the resulting murmuration functions. Numerically, the same phenomenon happens without the $N$ squarefree restriction (the coprime to $p$ restriction is asymptotically negligible), but the trace formula is more complicated and the murmuration conjecture has not been proved in this case.  

For the case of elliptic curves, again one numerically sees a sort of scale invariance in $\frac pX$, but there is also noise in the graph that does not seem to disappear, and a more reasonable conjecture is that the averages $\mathcal A_{\mathcal E}^{\pm}(p; X)$ only tend to a function in $\frac pX$ after  suitable smoothing (e.g., smooth by averaging over nearby $p$ for each $X$).

The other thing to notice about the graphs in \cref{fig:fixedroot1,fig:fixedroot2} is that the blue and red (i.e., $+$ and $-$) graphs are almost reflections of each other.  Indeed the limiting murmuration functions for opposite root numbers are precisely negatives of each other, so for asymptotic purposes one may simply consider the difference
\[  \mathcal A^\Delta_{\mathcal F}(p,X) = \mathcal A_{\mathcal F}^{+}(p,X) - \mathcal A_{\mathcal F}^{-}(p,X). \]
This has the advantage that the inner sum in the numerator for $\mathcal A^\Delta_{\mathcal F}(p,X)$ may now (in the case of $\mathcal F = \mathcal H_k$) be expressed as the trace
\[ \tr_{\mathcal F(N)} W T_p := \sum_{f \in \mathcal F(N)} w_f a_p(f) = \sum_{f \in \mathcal F^+(N)} a_p(f) -
\sum_{f \in \mathcal F^-(N)} a_p(f), \]
where $W$ is the Fricke involution on the newspace $S_k^\new(N)$ and $w_f$ is the root number of $f$.\footnote{Technically one usually considers the averages of $w_f a_p(f)$ over $\mathcal F$, rather than averaging separately over each root number, but this is approximately the same as $\frac 12\mathcal A^\Delta_F(p,X)$.}  Hence one can analyze such averages now with the trace formula.  We note that in the simple case of $k=2$, $N > 1$ squarefree coprime to $p$, the trace formula reads
\begin{equation} \label{eq:trWTp}
 \tr_{S_2^\new(N)} W T_p = \frac 12 \sum_{s^2 \le \frac{4p}{N}} H(s^2 N^2 - 4Np) - (p+1),
\end{equation}
where $s \in \Z$, and $H$ is the Hurwitz class number.

Here are some additional remarks: (1) This correlation between $a_p$'s and root numbers is something one only sees averaging over large families, and is quantitatively quite small --- after the normalization by $p^{1- \frac k2}$, each individual Fourier coefficient has size on the order of $\sqrt p$.  (2) The Birch and Swinnerton-Dyer conjecture asserts a subtle correlation between sizes of $a_p$'s and ranks of elliptic curves.  (3) For small $p$, a correlation between $a_p$ and the root number of newforms was already observed in \cite{me:pharis}; the ``limiting case'' of $p=1$ corresponds to the bias of root numbers toward $+1$ from \cite{me:refdim,me:rootno}. (4) One can more generally consider correlation of $a_n$'s, but for simplicity we stick to $a_p$'s --- this restriction makes patterns more apparent with less data.

\section{No root numbers} \label{sec:no-root}

One of the first variations that might come to mind is simply to consider averages without fixing root numbers, which is essentially the average of
$\mathcal A_{\mathcal F}^{+}(p,X)$ and $\mathcal A_{\mathcal F}^{-}(p,X)$.  As remarked above, they will tend to 0, but one can still look for structure in the behavior as $X \to \infty$, which can be thought of as a second-order term in comparing the convergence to murmuration functions for root number $+1$ versus root number $-1$.

In \cite{me:Wqbias}, I conjectured the existence such murmurations without root numbers for $\mathcal F = \mathcal H_k$, but not $\mathcal F = \mathcal E$.  Specifically, consider the weighted averages
\[ \mathcal A_{\mathcal F}(p, X) = p^{1 - \frac k2}
\frac{\sum'_{X \le N \le \beta X} \sqrt N \sum_{f \in \mathcal F(N)}  a_p(f)}{\sum'_{X \le N \le \beta X} \# \mathcal F(N)}, \]
where $\mathcal F(N)$ is the set of $f \in \mathcal F$ of level $N$.  Again, for simplicity, the prime on the outer sums means in our calculations we restrict to squarefree $N$ coprime to $p$, and we take $\beta = 2$.  Here we have inserted a scaling factor of $\sqrt N$ in the numerator to prevent the averages from going to 0.
 
\begin{figure}
\begin{minipage}{.5\textwidth}
\captionof{figure}{Weight 2 murmurations without root number for $X=2000$}
 \includegraphics[width=\textwidth]{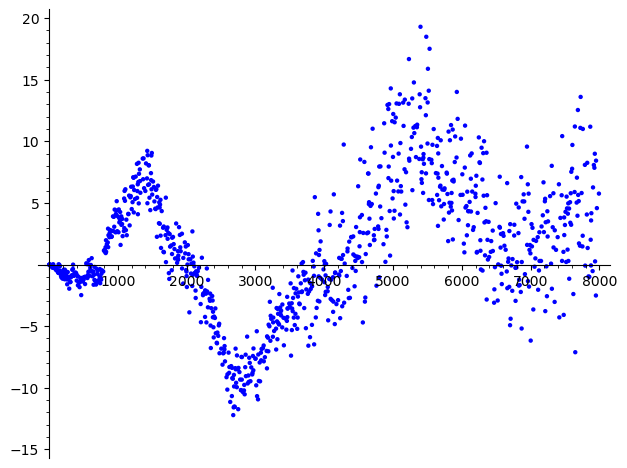}
\label{fig:no-root1}
\end{minipage}%
\begin{minipage}{.5\textwidth}
\captionof{figure}{Weight 2 murmurations without root number for $X=4000$}
 \includegraphics[width=\textwidth]{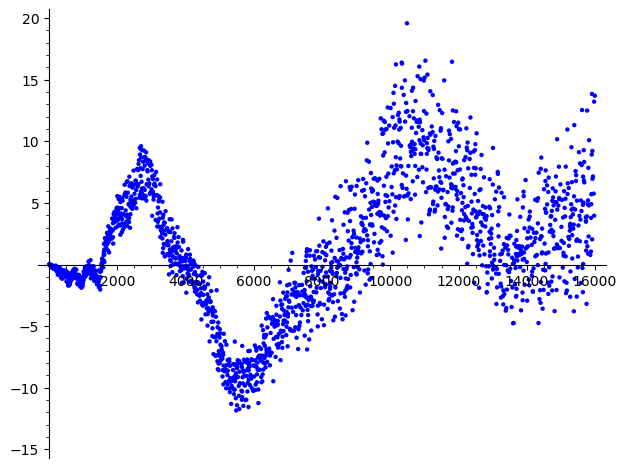}
\label{fig:no-root2}
\end{minipage}
\end{figure}

See \cref{fig:no-root1,fig:no-root2} for graphs of $\mathcal A_{\mathcal F}(p, X)$ with $\mathcal F = \mathcal H_2$, $p \le 4X$, where $X = 2000$ and $X= 4000$, respectively.  Note as in the original case of murmurations, there appears to be scale invariance in $\frac pX$, but there is more noise in the graphs.  Graphs without the $\sqrt N$ scaling look similar in shape, but have vertical range that tends to 0 as $X \to \infty$.   Graphs in higher weight look fairly similar, and for $\mathcal F = \mathcal H_k$ I conjectured that such graphs tend to a limiting murmuration function after appropriate smoothing.

\begin{figure}
\begin{minipage}{.5\textwidth}
\captionof{figure}{Elliptic curve $\sqrt N a_p$ averages without root number for $X=2000$}
 \includegraphics[width=\textwidth]{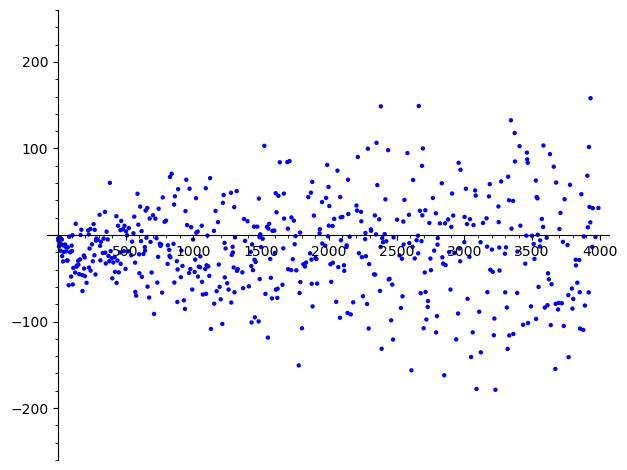}
\label{fig:Enoroot1}
\end{minipage}%
\begin{minipage}{.5\textwidth}
\captionof{figure}{Elliptic curve $\sqrt N a_p$ averages without root number for $X=4000$}
 \includegraphics[width=\textwidth]{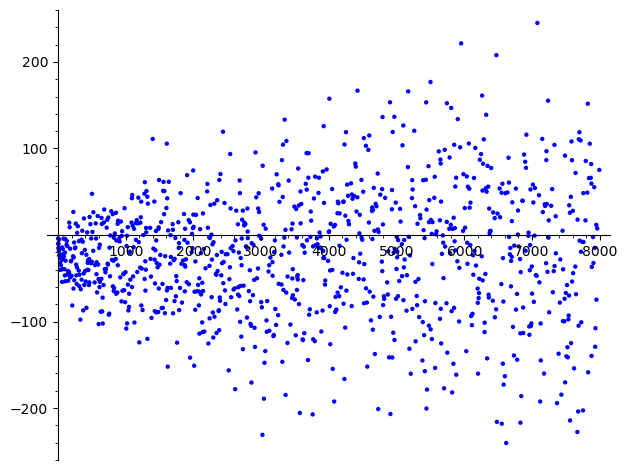}
\label{fig:Enoroot2}
\end{minipage}
\end{figure}

On the other hand, for elliptic curves, previous calculations of Sutherland indicate no such murmurations if one omits root numbers.  See \cref{fig:Enoroot1,fig:Enoroot2} for the elliptic curve analogues of \cref{fig:no-root1,fig:no-root2}, but with ranges restricted to $p \le 2X$.

To indicate the difference with the usual murmurations from \cref{sec:review}, we write down the trace formula for the inner sum $\tr_{\mathcal F(N)} T_p = \sum_{f \in \mathcal F(N)} a_p(f)$ in the simple case that $\mathcal F = \mathcal H_2$ and $N > 1$ is squarefree and coprime to $p$:
\begin{equation} \label{eq:trTp}
 \tr_{S_2^\new(N)} T_p = -\frac 12 \sum_{s^2 \le 4 p} \xi_{s^2- 4 p}(N) H(s^2 - 4p) + \mu(N) (p+1),
\end{equation}
where $\xi_\Delta$ is a certain multiplicative function defined in \cite{me:Wqbias}, and $\mu$ is the M\"obius function.  (If $q$ is a prime such that $q^2 \nmid \Delta$, then $\xi_\Delta(q) = {\Delta \leg q} - 1$.)

A significance difference is that in \eqref{eq:trTp} the number of class number terms grows with $p$, so is unbounded in a limit $p, N \to \infty$, whereas as in \eqref{eq:trWTp} the number of terms is bounded in $\frac pN$.  This makes  a theoretical analysis more challenging.

\section{M\"obius sums}
To try to analyze the averages of $a_p$'s without root numbers from \cref{sec:no-root} in the case of weight 2 modular forms, one might first try to analyze the contribution of the $\mu(N)(p+1)$ term in \eqref{eq:trTp}.  We remark that this term comes from removing the Eisenstein contribution from the trace of $T_p$ on $M_k(N)$.  Specifically, the contribution to $\mathcal A_{\mathcal H_2}(p,X)$ is
\begin{equation} \label{eq:mu-contrib} \mathcal A_{\mathcal H_2}^\mu(p,X) = (p+1) 
\frac{\sum'_{X \le N \le \beta X} \sqrt N \mu(N)}
{\sum'_{X \le N \le \beta X} \dim S_2^\new(N)}.
\end{equation}
For $N$ squarefree, $\dim S_2^\new(N) = \frac{\phi(N)}{12} + O(\log N)$, the approximate growth is
\begin{equation} \label{eq:mu-approx} \mathcal A_{\mathcal H_2}^\mu(p,X) \approx
\frac p{X^2} {\sideset{}{'}\sum_{X \le N \le \beta X} \sqrt N \mu(N)}.
\end{equation}
(Here $\approx$ means asymptotic up to a scalar.)

The right hand side is a modified version of Mertens function $M(X) = \sum_{1 \le N \le X} \mu(N)$ (or rather a weighted analogue of $M(\beta X) - M(X)$).  It is conjectured that $\lvert M(X)/\sqrt X \rvert$ is unbounded.  (This is subtle---recall that the Riemann hypothesis is equivalent to $M(X) = O(X^{\frac 12 + \epsilon})$.)  Similarly, one would expect that the sum $\sum' \sqrt N \mu(N)$ in the right hand side of \eqref{eq:mu-approx} should \emph{not} be $O(X)$, and which would mean $\mathcal A^\mu_{\mathcal H_2}(p,X)$ is \emph{not} bounded as $p, X \to \infty$ such that $\frac pX$ tends to a non-zero limit.

On the other hand, we just conjectured that $\mathcal A_{\mathcal F}(p, X)$ is bounded!  How can these be compatible?  Admittedly, the asymptotics of Mertens-like functions are quite subtle, so perhaps more care should be taken in the above approximation.  Still, we suggest that the erratic (in $X$) behavior of the M\"obius contribution $\mathcal A_{\mathcal H_2}^\mu(p,X)$ gets canceled out with the class number sums.  That is, we propose there is an almost magical interaction between the $\mu(N)$ term and the class number sum in $\tr T_p$ in \eqref{eq:trTp}.  Note that in $\tr T_p$ there are roughly $\sqrt p$ class numbers appearing, each of order approximately $\sqrt p$, and the $\mu(N)$ term is of size $p+1$.

\begin{figure}
\begin{minipage}{.5\textwidth}
\captionof{figure}{Class number versus M\"obius sums for $X=2000$}
 \includegraphics[width=\textwidth]{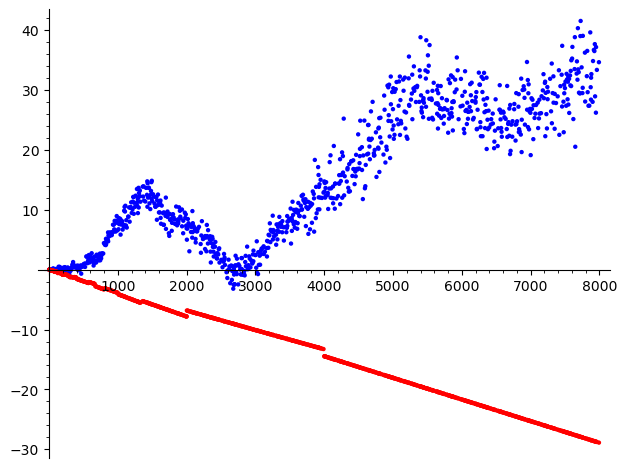}
\label{fig:mu1}
\end{minipage}%
\begin{minipage}{.5\textwidth}
\captionof{figure}{Class number versus M\"obius sums for $X=4000$}
 \includegraphics[width=\textwidth]{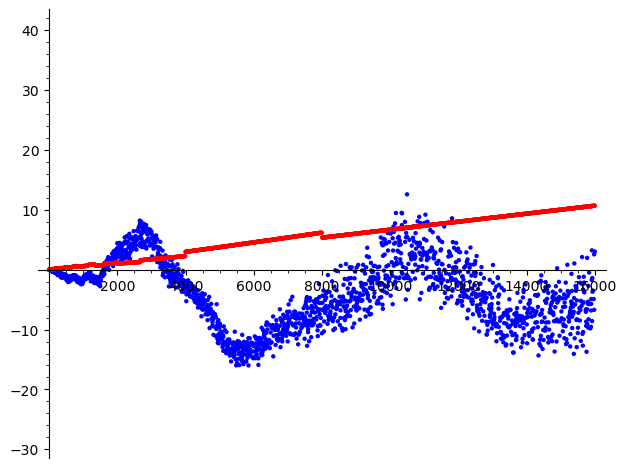}
\label{fig:mu2}
\end{minipage}
\end{figure}

\begin{figure}
\begin{minipage}{.5\textwidth}
\captionof{figure}{M\"obius contribution slopes}
 \includegraphics[width=\textwidth]{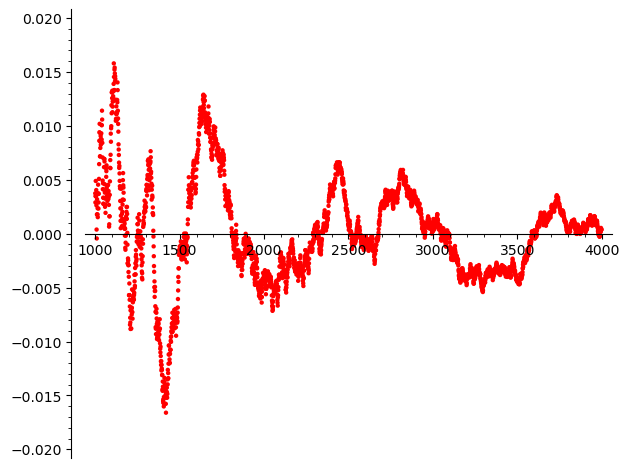}
\label{fig:muslope}
\end{minipage}%
\begin{minipage}{.5\textwidth}
\captionof{figure}{Murmurations on AL eigenspaces on $S_2(2q)$}
 \includegraphics[width=\textwidth]{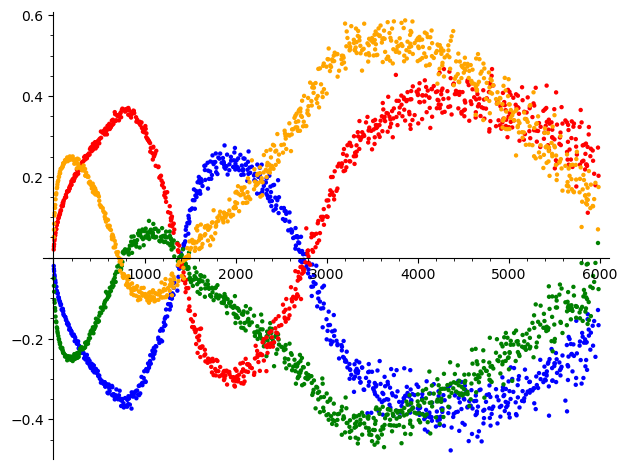}
\label{fig:pq-prof}
\end{minipage}
\end{figure}

As evidence for this interaction, we decompose the murmuration graphs in \cref{fig:no-root1,fig:no-root2} into the class number contribution (in blue) and the M\"obius contribution $\mathcal A_{\mathcal H_2}^\mu(p,X)$ (in red) in \cref{fig:mu1,fig:mu2}.  Namely, summing the red and blue graphs in the latter figures gives the scale-invariant murmration graphs in the former figures.  However, neither the class number nor the M\"obius contributions individually appear to be scale invariant, as the slope of the M\"obius contribution changes with $X$.  See \cref{fig:muslope} for a plot of the slopes of the M\"obius contribution for $1000 \le X \le 4000$.

Lastly, we remark upon the jaggedness of the red lines in these figures: note that the expression for the M\"obius contribution in \eqref{eq:mu-contrib} is not exactly linear in $p$ due to the fact that we restricted our sums over $N$ to $\gcd(N,p) = 1$.

\section{Local root numbers} \label{sec:loc-root}
Both the usual murmurations with respect to root numbers in \cref{sec:review} and without root number in \cref{sec:no-root} fall into a more general framework given in \cite{me:Wqbias}.  Namely, the averages considered in those two situations can be seen as weighted averages of quantities $\tr W T_p$ and $\tr T_p = \tr W_1 T_p$, respectively.

More generally one can consider averages of $\tr_{S_k^\new(N)} W_M T_p$'s, where $W_M$ is a suitably chosen Atkin--Lehner operator on $S_k^\new(N)$.  (In the averages, $M$ will potentially vary with $N$; in our notation, the Fricke involution is $W = (-1)^{k/2} W_N$.)  Conceptually, one can think of this variant as studying the correlation of Fourier coefficients with (the product of) local root numbers at some subset $S$ of places.  Looking at the global root number is the case where $S$ is all places of $\Q$, and no root number is the case where $S$ is the empty set.

We refer to \cite{me:Wqbias} for the full context and details, along with some theoretical results.  Here we just illustrate some examples in the simple case of levels of the form $N=pq$, where $p < q$ are distinct primes.  We propose that a good way to study murmurations with respect to local root numbers is to look at murmuration graphs for each Atkin--Lehner eigenspace, analogous to our original murmuration pictures
in \cref{fig:fixedroot1,fig:fixedroot2} where we plotted graphs for each global root number separately.

For a newform $f$ of level $N=pq$, the global root number $w_f = \prod_v w_{f,v}$ where $w_{f,v}$ is the local root number at $v$.  Here $w_{f,\infty} = (-1)^{k/2}$,  $w_{f,p}, w_{f,q} \in \{ \pm 1 \}$, and $w_{f,v} = +1$ for all other $v$.  There are 4 Atkin--Lehner eigenspaces, which we denote by the 4 sign patterns \verb|++|,  \verb|+-| and \verb|-+|, \verb|--|.  E.g., \verb|+-| will refer to the subspace of $S_k^\new(N)$ generate by newforms with $w_{f,p} = +1$ and $w_{f,q} = -1$.

\begin{figure}
\begin{minipage}{.5\textwidth}
\captionof{figure}{Murmurations on AL eigenspaces on $S_4(2q)$}
 \includegraphics[width=\textwidth]{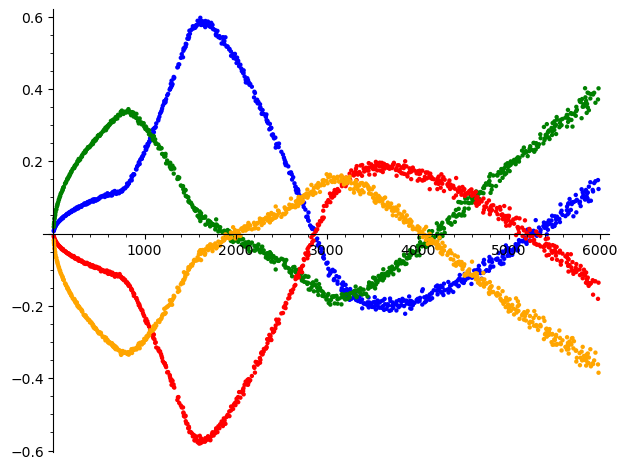}
\label{fig:pq-prof4}
\end{minipage}%
\begin{minipage}{.5\textwidth}
\captionof{figure}{Murmurations on AL eigenspaces on $S_4(pq)$}
 \includegraphics[width=\textwidth]{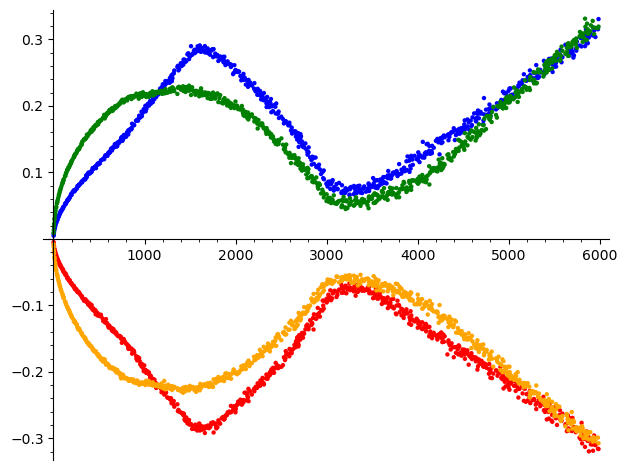}
\label{fig:pq-all-prof4}
\end{minipage}
\end{figure}

We present two types of examples: (i) we fix $p$ and vary $q$, and (ii) we vary both $p$ and $q$. 
See \cref{fig:pq-prof,fig:pq-prof4} for a graphs of averages of $a_\ell$'s over each Atkin--Lehner eigenspace in $S^\new_2(pq)$ and $S^\new_4(pq)$.  Here $p = 2$ is fixed and $3000 < q < 6000$.    See \cref{fig:pq-all-prof4} for $S^\new_4(pq)$ where $p < q$ and $6000 < pq < 12000$.

In all of these plots, the blue and green dots correspond to signs \verb|++| and \verb|--| and red and orange dots to signs \verb|+-| and \verb|-+|, respectively. 
Adding all 4 colors corresponds to looking at murmurations with no root number.  Adding just the blue and green (resp., red and orange) graphs corresponds to looking at murmurations for global root number $+1$ (resp., $-1$) when $k=4$, and $-1$ (resp.\ $+1$) when $k=2$.  
 
We remark that if one goes beyond squarefree levels, one can also look at something more refined than just the local root number at ramified places --- one can also look at the local inertial type of the representation.  Trace formulas that will allow us to do such calculations have been worked out in some cases recently in \cite{knightly}.
 
\section{Class number sums}
Zubrilina \cite{zubrilina} proved murmurations for modular forms by deriving asymptotics for short class number sums of the form
\[ \sideset{}{'}\sum_{X \le N \le Y} H(s^2 N^2 - 4Np). \]
This is a sum over suitable levels for a given $s$-term in \eqref{eq:trWTp}.  One can view this as a sum of class numbers $H(f_s(N;p))$ where $f_s$ is a polynomial, and $N$ runs over a restricted range.  The existence of murmurations amounts to suitably normalized sums of class numbers having a limit as $p, N \to \infty$ such that $\frac pN \to x$ for any $x > 0$.

Numerically we have observed one has similar behavior if one slightly varies the polynomial $f_s(N;p)$, e.g., by varying scalars or adding or subtracting smaller order terms.  We do not have a comprehensive philosophy of how such class number sums should behave, or what polynomials are interesting, so we simply present a few numerical examples.  We begin by restricting to quadratic polynomials $f_s(N;p)$.

\begin{figure}
\begin{minipage}{.33\textwidth}
\captionof{figure}{Class number sum for $s^2 N^2 - 4Np$}
 \includegraphics[width=\textwidth]{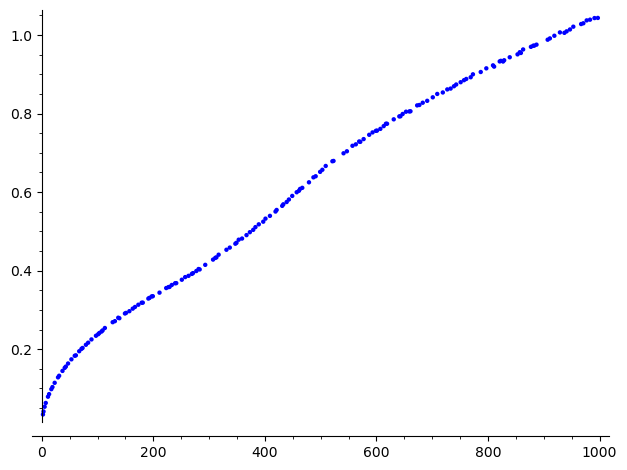}
\label{fig:cl0avg}
\end{minipage}%
\begin{minipage}{.33\textwidth}
\captionof{figure}{Class number sum for $s^2 p^2 - 4Np$}
 \includegraphics[width=\textwidth]{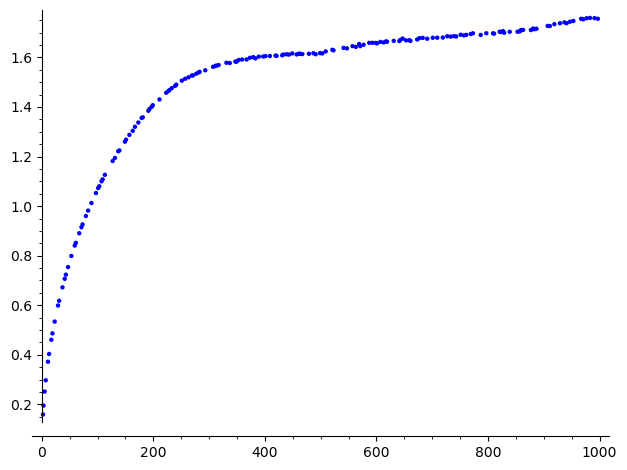}
\label{fig:cl1avg}
\end{minipage}%
\begin{minipage}{.33\textwidth}
\captionof{figure}{Class number sum for $s^2 + 1 - 3Np$}
 \includegraphics[width=\textwidth]{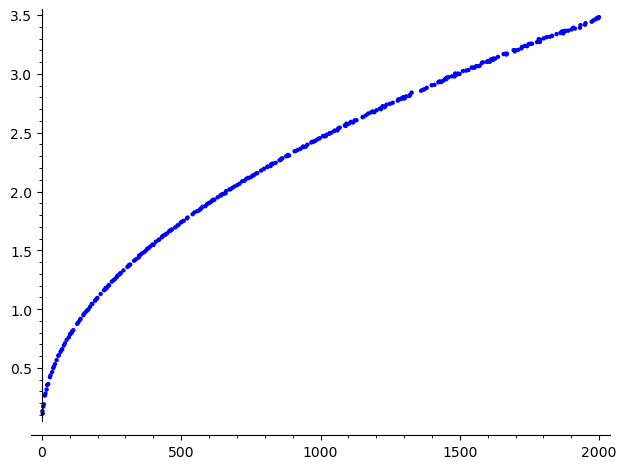}
\label{fig:cl2avg}
\end{minipage}%
\end{figure}

As a baseline, in \cref{fig:cl0avg}, we present a graph of
\[ A_0(p,X) := \frac{\sum'_{X \le N \le 2X} \sum_{s=0}^1 H(s^2 N^2 - 4 Np)}{\sum'_{X \le N \le 2X} N}, \]
for $X=1000$ and $p \le X$.  Here $\sum'$ denotes a restriction to squarefree $N$.  This is essentially the $s=0$ and $s=1$ term contributions to murmurations for $\mathcal F = \mathcal H_2$ from \eqref{eq:trWTp}.  We chose to include both $s=0$ and $s=1$ terms so that one can see some oscillation --- the terms for larger $s$ do not contribute in the range $p \le X$.

In \cref{fig:cl1avg}, we simply consider $f_s(N;p) = s^2 p^2 - 4 Np$ instead of $s^2 N^2 - 4Np$.  Here larger $s$ terms come into play in the range $p \le X$, and we plot
\[ A_1(p,X) := \frac{\sum'_{X \le N \le 2X} \sum_{s=0}^4 H(s^2 p^2 - 4 Np)}{\sum'_{X \le N \le 2X} N} \]
for $X = 1000$.  Again there is some oscillatory behavior.

In \cref{fig:cl2avg}, we plot
\[ A_2(p,X) := \frac{\sum'_{X \le N \le 2X} \sum_{s=0}^4 H(s^2 + 1 - 3 Np)}{\sum'_{X \le N \le 2X} N} \]
for $X=1000$.  Here each class number discriminant is about size $-3Np$, and it seems there is not enough variation in adding $s^2+1$ to cause  oscillation.

In all of the above examples, while we only included graphs for $X=1000$, graphs look similar for other values of $X$, i.e., further calculations indicate  scale invariance in $\frac pX$, i.e., there is a limiting graph as $X \to \infty$.

\begin{figure}
\begin{minipage}{.5\textwidth}
\captionof{figure}{Class number sum for $s^2 N^3 - 4N^2p$}
 \includegraphics[width=\textwidth]{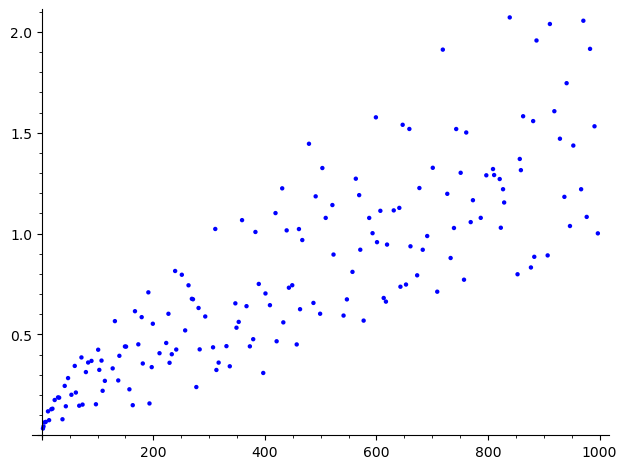}
\label{fig:cl3avg}
\end{minipage}%
\begin{minipage}{.5\textwidth}
\captionof{figure}{Class number sum for $s^2N^3 -tN- 4N^2p$}
 \includegraphics[width=\textwidth]{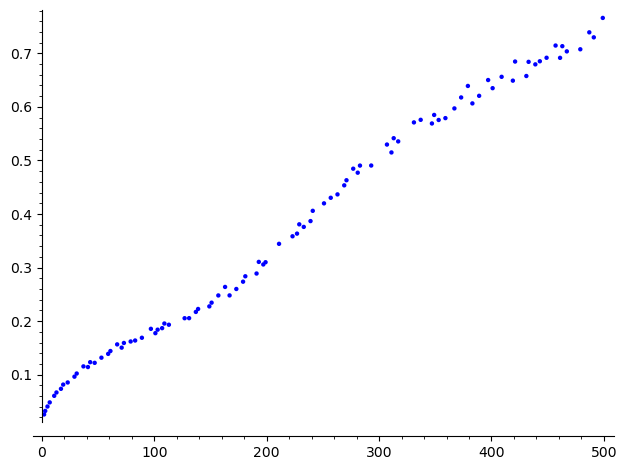}
\label{fig:cl4avg}
\end{minipage}%
\end{figure}

Now we come to an example of a cubic polynomial $f_s(N;p)$.  We plot
\[ A_3(p,X) := \frac{\sum'_{X \le N \le 2X} \sum_{s=0}^1 H(s N^3 - 4 N^2p)}{\sum'_{X \le N \le 2X} N^{3/2}} \]
for $X=1000$ in \cref{fig:cl3avg}.
Note that we changed the normalization factor in the denominator to account for the class numbers now being of size approximately $N^{3/2}$ (assuming $p \approx X$).  Again the graphs look roughly similar as $X \to \infty$, but there does not appear to be convergence to an actual function, because we are averaging too sparse a set of class numbers.

But we can modify the last example by averaging over more class numbers.  For instance, instead of average over $\approx X$ class numbers for discriminants of size $\approx X^3$, we consider the following average over $\approx X^{3/2}$ such class numbers,
\[ A_4(p,X) := \frac{\sum'_{X \le N \le 2X} \sum_{t \le \sqrt X}\sum_{s=0}^1 H(s^2 N^3 - tN - 4 N^2 p)}{\sum'_{X \le N \le 2X} \sum_{t \le \sqrt X} N^{3/2}}. \]
This average seems to converge to a function in $\frac pX$ --- see \cref{fig:cl4avg} for a plot of $A_4(p,X)$ with $X = 500$.

\section{$L$-values} \label{sec:Lval}

Several of the murmurations we have considered are just averages of the geometric side of a trace formula.  There are other kinds of trace formulas one can consider: the Kunzetsov trace formula, the Petersson trace formula, and more generally relative trace formulas.  The Kuznetsov trace formula will give sums of Fourier coefficients weighted by inverse Petersson norms.  We will not address that, but just discuss one different relative trace formula.

One of Jacquet's first relative trace formulas relates toric periods on quaternion algebras to twisted central values of $L$-functions, reproving a result of Waldspurger.  Evaluating the geometric side for the quaternion algebra gives an exact formula for average $L$-values 
\[ \Lambda(N,k,D,n) = \sum_{f \in \mathcal H_k(N)} \frac{L(1/2,f)L(1/2,f \otimes \chi_{-D})}{(f,f)} a_n(f), \]
where $\chi_{-D} = {-D \leg \cdot}$, for some fixed fundamental discriminant $-D < 0$.  This was carried out in \cite{FW,me:doub-avg} under some hypotheses on $N$ and $D$.

We just indicate what the geometric side looks like in the following simple setting, which was originally considered in \cite{MR} (using the Gross--Zagier formula, rather than Jacquet's relative trace formula).  Suppose $0 < D \equiv 3 \mod 4$, and $N \ne p$ is a prime inert in $K = \Q(\sqrt {-D})$.  Let $h_{D} = h_K$ and $u_{D} = [\calO_K^\times : \Z^\times]$.  Also assume $k=2$ for simplicity.  Then
\begin{equation} \label{eq:MR} \frac{\sqrt D u_D^2}{2 \pi} \Lambda(N,k,D,n) = \frac{12 h_{D}^2}{N-1} \sigma_N(n) + u_D r(nD)h_{D} + u_D^2 \sum_{m=1}^{\lfloor nD/N \rfloor} \Phi(m,N),
\end{equation}
where $\sigma_N(n)$ is the sum of divisors $d \mid n$ such that $(d,N) = 1$, $r(nD)$ is the number of ideals of norm $nD$ in $\calO_K$ and 
\[ \Phi(m,N) = d((m,D)) r(m) r(pD - mN). \]  
(Here $d = \sigma_0$ is the number-of-divisors function.)

There are a couple of ways one could investigate murmuration analogues for these quantities.  For the standard trace formula for $\tr T_p$ or $\tr W_N T_p$ on $S_k^\new(N)$, there are 3 parameters one can vary: $N, k$ and $p$, and we varied $p$ and $N$.  (One can also consider varying $k$---see
\cite{BBLLD}.)  Here there are 4 parameters: $N, k, D, n$.  We choose to keep $k$ fixed and average over $N$ as in the original murmurations setting, and vary either $D$ or $n=p$.  

\begin{figure}
\begin{minipage}{.5\textwidth}
\captionof{figure}{Averages of twisted $L$-values for $X = 4000$}
 \includegraphics[width=\textwidth]{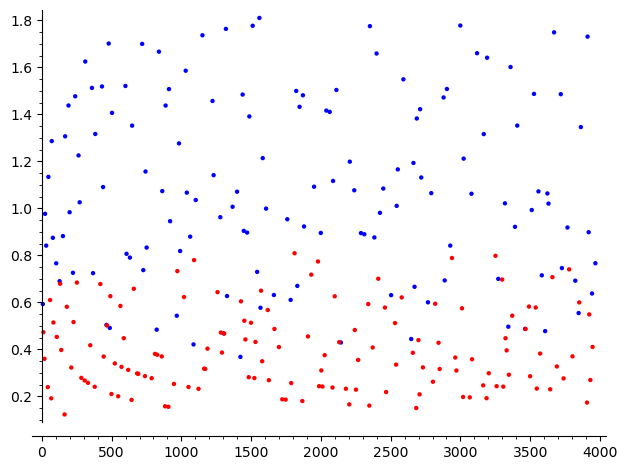}
\label{fig:LvalD}
\end{minipage}%
\begin{minipage}{.5\textwidth}
\captionof{figure}{Averages of twisted $L$-values for $X = 4000$}
 \includegraphics[width=\textwidth]{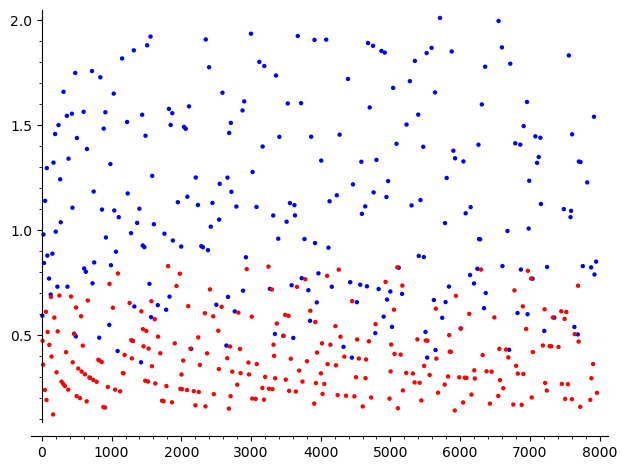}
\label{fig:LvalD2}
\end{minipage}
\end{figure}

First consider varying $D$.  This amounts to considering murmurations for squares of Fourier coefficients of half-integral weight forms (see \cite{wald81}).  
In \cref{fig:LvalD,fig:LvalD2}, we plot averages 
\[ A_{L,D}(D,X) = \frac{\sum'_{X \le N \le 2X} \Lambda(N,2,D,1)}
{\sum'_{X \le N \le 2X} 1} = \frac{\sum'_{X \le N \le 2X}\sum_{f \in \mathcal H_2(N)} \frac{L(1/2,f)L(1/2,f \otimes \chi_{-D})}{(f,f)} }
{\sum'_{X \le N \le 2X} 1}  \]
for $X = 4000, 8000$, where $-D \equiv 1 \mod 4$ is a negative prime fundamental discriminant.  Here the prime in the sum over $N$ refers to a restriction to primes inert in $K = \Q(\sqrt{-D})$ so that \eqref{eq:MR} is valid.
We take $3 < D \le X$ along the horizontal axis, and color the $-D \equiv 1 \mod 8$ points blue and the $-D \equiv 5 \mod 8$ points red.  

We tried a few different normalizations/weightings of these averages, and while the graphs look roughly similar as $X$ increases, we did not find a simple one which is scale invariant in $\frac DX$, or is clearly indicative of murmurations.  Note that we divide by the number of levels in the average, rather than the number of newforms, to prevent the graphs from shrinking in scale as $X$ grows.  We also remark that such averages only pick up forms with root number $+1$ because of the $L(1/2,f)$ factor.

\begin{figure}
\begin{minipage}{.5\textwidth}
\captionof{figure}{Averages $a_p$'s weighted by $L$-values for $X=2000$}
 \includegraphics[width=\textwidth]{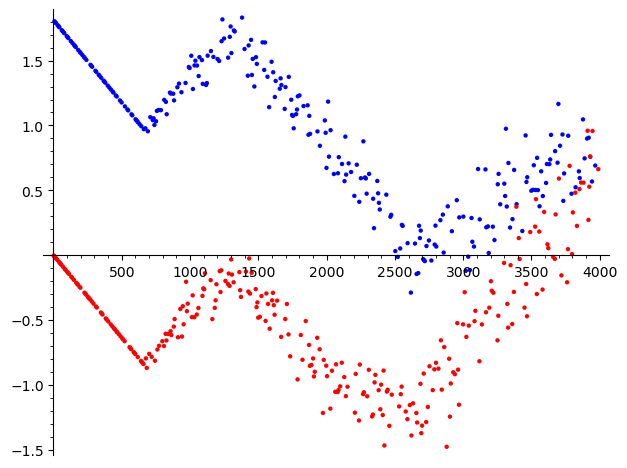}
\label{fig:Lval}
\end{minipage}%
\begin{minipage}{.5\textwidth}
\captionof{figure}{Averages $a_p$'s weighted by $L$-values for $X=4000$}
 \includegraphics[width=\textwidth]{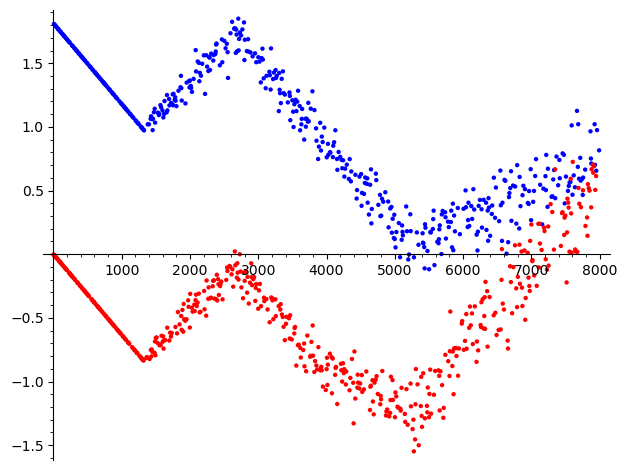}
\label{fig:Lval2}
\end{minipage}
\end{figure}

Finally, we consider the case of fixing $D$ and varying $n=p$.  Here we consider averaging the scaled quantity in \eqref{eq:MR}.  Let
\[ A_{L,p}(p,X) = \frac{\sum'_{X \le N \le 2X} \Lambda(N,2,3,p)}
{\sum'_{X \le N \le 2X} 1} = \frac{\sum'_{X \le N \le 2X}\sum_{f \in \mathcal H_k(N)} \frac{L(1/2,f)L(1/2,f \otimes \chi_{-3})}{(f,f)} a_p }
{\sum'_{X \le N \le 2X} 1}.  \]
In \cref{fig:Lval,fig:Lval2}, we plot $A_{L,p}(p,X)$ for $p \le X$ ($p \ne 3$) and $X = 2000, 4000$, coloring the values blue or red according to $\chi_{-3}(p) = +1$ or $-1$.  

These graphs are apparently scale invariant in $\frac pX$, and we expect the blue and red graphs to tend to a limiting function after appropriate smoothing.  The first linear part of the graph corresponds to the so-called \emph{stable range} $p < \frac ND$ where the sum on the right hand side of \eqref{eq:MR} vanishes (in this range, it is not hard to prove the limit exists).  This is somewhat analogous to murmurations in \cref{sec:review}, where the first part of the murmuration graph is simply given by a square root function (or a linear function if one graphs in terms of $\sqrt{\frac pX}$, which turns out to be nicer), and more trace formula terms contribute the further to the right one goes on the graph.  

The reason to separate the cases of $p$ inert or split in $K = \Q(\sqrt{-D})$ it because that affects the second term on the right hand side of \eqref{eq:MR}.
We also note that $\sum' 1$ is the correct scaling in the denominator to get scale invariance in $\frac pX$ because the main term (the first term on the right) of \eqref{eq:MR} is of size $\approx \frac pN$.

\section{Representations by quadratic forms}
For an elliptic curve $E$ at a good prime $p$, $a_p(E) = p+1 - \# E(\F_p)$, which is the deviation of expected minus actual number of solutions mod $p$.  Hence murmurations for elliptic curves indicate how the error in first-order point count estimates mod $p$ is correlated with the conductor.

One other situation in number theory where the error of first-order solution count estimates is well studied is representations by quadratic forms.  In fact, both traces of $T_p$'s and the averages $A_{L,p}$ of $L$-value weighted $a_p$'s can be interpreted in terms of Brandt matrix entries, which can be expressed as representations numbers of quadratic forms (e.g., see \cite[Section 3.3]{me:doub-avg}).

The most classical case is that of representations of integers by binary quadratic forms, which includes the Gauss circle problem.  Let $-D < 0$ be a discriminant, and $r_D(n)$ be the total number of ways to represent $n$ by a positive definite reduced binary quadratic form of discriminant $-D$.  Then
\[ r_D(n) = 2u_D \sum_{d \mid n} \chi_{-D}(n). \]
In particular, for $n=p$ the $\chi_{-D}(p)$ measures the error of the first-order approximation for the number of ways to represent $p$ by a reduced form of discriminant $-D$.  
Thus one analogue of the murmuration averages for elliptic curves in the setting of positive definite binary quadratic forms is the quantity
\[ A^{\mathcal D}_{\textrm{BQF}}(p,X) = \frac 1{\sqrt X} \sideset{}{'}\sum_{\substack{X \le D \le 2X \\ -D \in \mathcal D }} \chi_{-D}(p),  \]
where $\mathcal D$ is a chosen class of fundamental discriminants (see below) and the sum is taken over $D$ coprime to $p$.  For a given $p$, the sum on the right should be roughly of size $\sqrt X$, which is why we normalize by $\frac 1{\sqrt X}$.

\begin{figure}
\begin{minipage}{.5\textwidth}
\captionof{figure}{Binary quadratic form murmurations with $D$ odd for $X=5000$}
 \includegraphics[width=\textwidth]{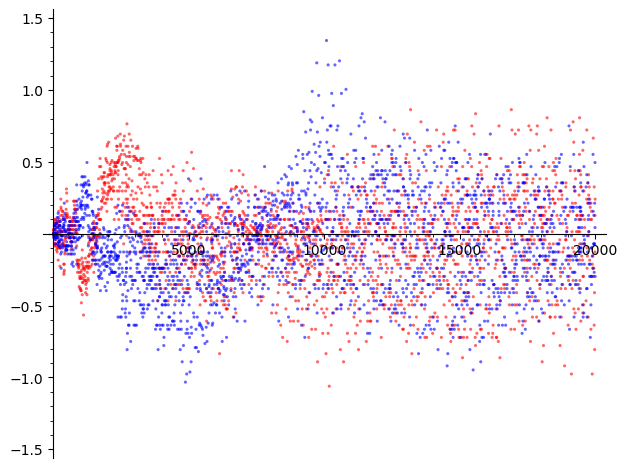}
\label{fig:bqf1odd}
\end{minipage}%
\begin{minipage}{.5\textwidth}
\captionof{figure}{Binary quadratic form murmurations with $D$ odd for $X=10000$}
 \includegraphics[width=\textwidth]{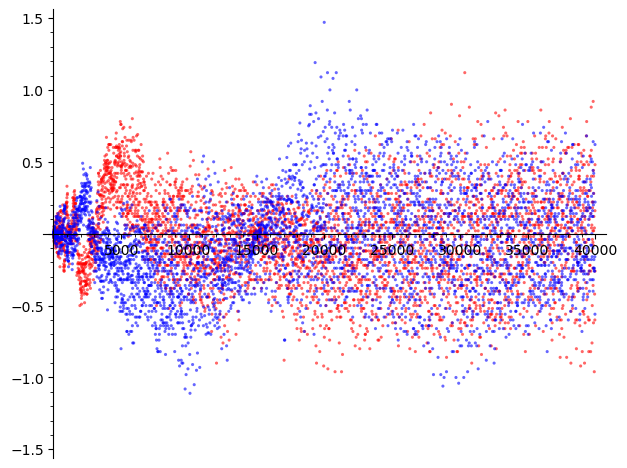}
\label{fig:bqf2odd}
\end{minipage}
\end{figure} 

\begin{figure}
\begin{minipage}{.5\textwidth}
\captionof{figure}{Binary quadratic form murmurations with $D$ even for $X=5000$}
 \includegraphics[width=\textwidth]{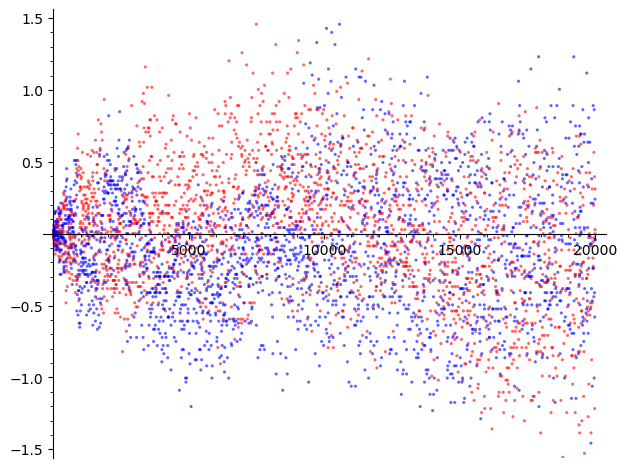}
\label{fig:bqf1even}
\end{minipage}%
\begin{minipage}{.5\textwidth}
\captionof{figure}{Binary quadratic form murmurations with $D$ even for $X=10000$}
 \includegraphics[width=\textwidth]{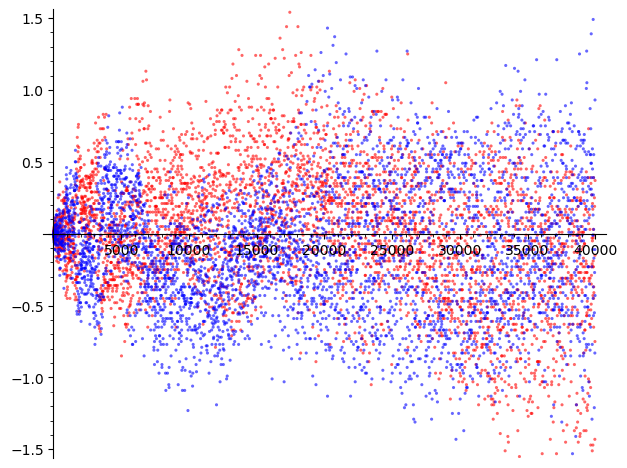}
\label{fig:bqf2even}
\end{minipage}
\end{figure} 

In \cref{fig:bqf1odd,fig:bqf2odd} we plot the quantities $A^{\mathcal D}_{\textrm{BQF}}(p,X)$ for $p \le 4X$, where $\mathcal D$ is either the set of fundamental discriminants $-D \equiv 1 \mod 8$ (in blue) or $-D \equiv 5 \mod 8$ (in red).  The first plot is with $X=5000$ and the second is with $X=10000$.  See \cref{fig:bqf1even,fig:bqf2even} for analogous plots where $\mathcal D$ is set of fundamental discriminants of the form $-4d$ with either $d \equiv 1 \mod 4$ (in blue) or $d \equiv 3 \mod 4$ (in red).
While there is a lot of noise in these graphs, they appear to be roughly scale invariant, which indicates the existence of murmurations in this setting.

Indeed, murmurations for quadratic (and general) Dirichlet characters were already established in \cite{LOP} under GRH.  They smooth out the noisiness by averaging over nearby primes.  Very recently, Cowan \cite{cowan:murm} obtained unconditional murmuration results, which also apply to the setting of quadratic Dirichlet characters.

%
%

\end{document}